\documentclass[12pt,reqno]{amsart}
\usepackage{amsmath,amsfonts,amssymb,amscd,amsthm,amsbsy,epsf}
\usepackage[dvips]{graphicx}
\usepackage{comment}
\usepackage[dvips]{color}
\usepackage{xcolor}

 \usepackage{lmodern}

\numberwithin{equation}{section}
\newtheorem{theorem}{Theorem}
\newtheorem{meta-thm}[theorem]{Meta-Theorem}
\newtheorem{lemma}[theorem]{Lemma}

\newtheorem{proposition}[theorem]{Proposition}

\newtheorem{definition}[theorem]{Definition}
\theoremstyle{remark}
\newtheorem{remark}[theorem]{Remark}

\def\POC{{Periodic Orbit Condition}}

\def\Id{\operatorname{Id}}
\def\Im{\operatorname{Im}}

\def\A{{\mathcal A}}
\def\cA{{\mathcal A}}
\def\B{{\mathcal B}}

\def\G{{\mathcal G}}
\def\cE{{\mathcal E}}
\def\E{{\mathcal E}}

\def\M{{\mathcal M}}

\def\complex{{\mathbb C}}
\def\integer{{\mathbb Z}}
\def\nat{{\mathbb N}}
\def\NN{{\mathbb N}}

\def\CC{{\mathbb C}}
\def\real{{\mathbb R}}
\def\RR{{\mathbb R}}
\def\torus{{\mathbb T}}

\def\ZZ{{\mathbb Z}}

\def\eps{\varepsilon}
\def\th{\theta}
\def\dt{\delta}
\def\si{\sigma}
\def\la{\lambda}
\def\ga{\gamma}
\def\halpha{ {\hat \alpha}}
\def\hbeta{ {\hat \beta}}

\def\<{\langle}
\def\>{\rangle}

\begin{document}
\title[Coboundaries over integrable systems]{Nonconmutative coboundary equations over integrable systems} 
\author[R. de la Llave]{Rafael de la Llave 
}
\address{
School of Mathematics,
Georgia Institute of Technology,
686 Cherry St. Atlanta GA. 30332-1160, USA }
\email{rafael.delallave@math.gatech.edu}
\thanks{R.L. has been partially supported by NSF grant DMS 1800241.}

\author[M. Saprykina]{Maria Saprykina 
}
\address{
  Department of Mathematics, KTH Royal Institute of Technology, Stockholm, Sweden}
\email{masha@kth.se}





\begin{abstract}
We prove an analog of Liv\v{s}ic theorem for real-analytic families of
cocycles over an integrable system with values in a Banach algebra
$\G$
or a Lie group.  

Namely, we consider an
integrable dynamical system $f:\M \equiv\torus^d \times [-1,1]^d\to
\M$, $f(\theta, I)=(\theta + I, I)$, and a real-analytic family of
cocycles $\eta_\eps : \M \to \G$, indexed by a complex parameter $\eps$ in an open ball $\cE_\rho \in\CC$. We show that if
$\eta_\eps$ has trivial periodic data, i.e., 
$$
\eta_\eps(f^{n-1}(p))\dots \eta_{\eps} (f(p))\cdot \eta_{\eps} (p)=Id
$$ 
for each periodic point $p=f^n p$
and each $\eps \in \cE_{\rho}$, then there exists a real-analytic family of
maps $\phi_\eps: \M \to \G$ satisfying the coboundary equation
$$
\eta_\eps(\theta, I)=\phi_\eps^{-1}\circ f(\theta, I)\cdot \phi_\eps
(\theta, I)
$$ 
for all $(\theta, I)\in \M$ and $\eps \in \cE_{\rho/2}$.

We also 
show that if the coboundary equation above with an analytic left-hand side $\eta_\eps$ 
has a solution in the sense of formal power series in $\eps$, then it has an analytic solution.

\end{abstract} 
\maketitle


\section{Introduction}\label{sec:intro}

On the annulus $\M =\torus^d \times [-1,1]^d$ consider an
integrable dynamical system
$f:\M \to\M$ of the form  $f(\theta, I)=(\theta + I, I)$
and a real-analytic family of 
maps $\eta_\eps :  \mathcal M \to \G$,
where $\eps\in \cE_\rho$ is a complex one-dimensional parameter, and $\G$ is a Banach algebra 
(with small modifications we can let $\G$ be a Lie group, see the last section). Assume that $\eta_0=\Id $.
We ask the question of whether $\eta_\eps$ is a {\it coboundary}, i.e., whether  there exists an analytic family 
$\phi_\eps: \mathcal M \rightarrow \G$,  $\eps \in \cE_{\rho/2}$,
solving the following coboundary equation:
\begin{equation} \label{coboundary} 
  \eta_\eps (\th,I) = \phi^{-1}_\eps \circ f(\th,I)  \cdot  \phi_\eps(\th, I),
\end{equation}
where the
dot between two elements of $\G$ denotes the
product in $\G$.

\bigskip

%
Equation  \eqref{coboundary}  can be expressed by saying  that cocycle $\eta_\eps$ is conjugated (equivalent)
to the identity, which coincides with $\eta_0$, and $\phi_\eps$ should be seen as the conjugacy.  
In this sense, the present paper is part of the rigidity program for cocycles. 

Coboundary  equations of type \eqref{coboundary} appear naturally in many problems in dynamics, in particular, in the study of rigidity for integrable Hamiltonian systems. 

\bigskip

An obvious necessary condition for $\eta_\eps$ 
to be a coboundary is the following \POC\  (POC for short). 
\begin{definition}
  \label{def:POC-general}
We say that $\eta_\eps$ satisfies the \POC\  (POC) if for every $\eps \in \E_\rho$ and for every 
$p \in \M$ such that $p = f^N(p)$ we have:
\begin{equation}\label{POC-general}
  \eta_\eps\circ f^{N-1}(p)\cdot \, \cdots\, \cdot    \eta_\eps \circ f (p) \cdot  \eta_\eps(p)  = \Id .
\end{equation}
\end{definition}

\bigskip

Our main result, Theorem \ref{main}, states that if this obvious necessary condition is
met, then there is an analytic solution of the coboundary
equation above. 

As a related result,   we
show that if the coboundary equation \eqref{coboundary} with an analytic left-hand side $\eta_\eps$ 
has a solution in the sense of $C^0$ formal  power series, then it has an analytic solution (which may be different from the original one). This provides a very dramatic
bootstrap of the regularity of the solutions to \eqref{coboundary}.

Both of the above statements remain valid when the family of cocycles
takes values in a Lie group. See Section~\ref{sec:Lie}.

\bigskip


The results of the present paper are related to the notion of the \emph{uniform integrability}  introduced by 
Poincar\'e in 
\cite[Chapter V]{Poincare99}.
 An analytic family of Hamiltonians
is called uniformly integrable if it can be reduced to an integrable
one by a canonical change of variables.  In the cited chapter, Poincar\'e described 
certain obstructions (i.e., necessary conditions) to the
existence of such families of  changes of variables even in the sense of formal
power series. These conditions can be expressed via certain integrals over periodic orbits. 
Poincar\'e  extended slightly these conditions and verified them in concrete problems such as the 3-body
problem. In \cite{Llave96} it is proved that these conditions are also sufficient: 
vanishing of the above obstructions to  integrability in the sense of
power series implies
the existence of an analytic integrating change of variables.



In the present paper we consider an analytic family of non-commutative coboundary equations \eqref{coboundary}, where the base dynamics, $f$, is an
integrable system (parabolic dynamics).
We prove the existence of an analytic family of solutions to this equation under the \POC\ \eqref{POC-general} (which is the obvious  necessary condition for the existence of such solutions).

As an iterative step of the proof we solve the linearized (commutative) coboundary equation, producing  solutions with tame estimates in the sense of Nash-Moser theory.  Recall that, in the context of analytic spaces,
an operator $\mathcal L$ is said to satisfy {\it tame estimates}
in a monotonic class of domains $D_\rho$, when there exists $\tau>0$ such that for
all $0 < \rho \le \rho_0$, given
$\psi: D_\rho \rightarrow \G$, we
have $\mathcal L \psi$ is analytic  in $D_{\rho - \delta}$ for
  every $0 < \delta < \rho$ and
  \[
  \|\mathcal L \psi \|_{\rho - \delta} \le C \delta^{-\tau} \|\psi\|_\rho .
  \]
(A similar notion of tameness exists in finite regularity
spaces.)
Informally speaking, our KAM iterative method shows
that we can pass from the tame estimates on the solutions of
commutative coboundary equations  to the analytic solutions of non-commutative
cobundary equations
for families (under the \POC).

\medskip

The type of statements where the Periodic Orbit Condition implies the existence of solutions to a functional equation is often referred to as Liv\v{s}ic theorems due to the seminal papers \cite{Livsic71, Livsic72}. 
There is a whole spectrum of contexts in which these results have been proved over the past decades: in particular, different classes of
smoothness, types of the base dynamics, types of the sufficient conditions, types of the coboundary equation (i.e., commutative or not) and others. Let us  mention just a few results that are closest to ours. 

\medskip

In the case when the base dynamics $f$
is {\it hyperbolic}, a lot is known. The most famous works are the classical ones by Livsic, \cite{Livsic71, Livsic72}, where the base dynamics is transitive Anosov, and the solutions are studied in low regularity. 
  Some generalizations and improved
results on the regularity were obtained in \cite{LlaveMM86}; in particular, they show that, in
appropriate spaces, there is no loss of regularity while solving the commutative coboundary equation.   
There are other methods 
  \cite{NiticaT95, NiticaT01, NiticaT02, KatokN11} that eliminate
  the use of power series and rely only on periodic orbit condition.

 In the analytic regularity, the solvability of the commutative equation is done in
 \cite{Llave97}, but the estimates presented there are not tame. Note that the notion of tame estimates depends very much on
the families of domains used.The estimates presented  in \cite{Llave97}
are not tame in the families of domains considered there, but it seems possible that they are tame in other families of domains.

The case of  non-commutative coboundary equations
over Anosov systems was clarified in
\cite{LlaveW10}. 

\medskip 

For commutative
cohomology equations over integrable systems ({\it parabolic} base dynamics) the regularity of solutions under
a variant of the \POC\ was studied in \cite{Llave96}.
In that paper it was shown that there are solutions satisfying tame estimates.

\medskip  

Another case when the cohomology equation has been studied is
the case of {\it quasi-hyperbolic} automorphisms of the torus. The paper \cite{Veech86}
established that, under the \POC, the solutions of
the equation with $C^\infty$ data are $C^\infty$. 


\medskip  
  
The assumptions of the present paper require that a whole analytic {\it family} 
satisfies the \POC.  Nevertheless, there are other methods: 
  \cite{LlaveMM86, LlaveW10,NiticaT95, NiticaT01, NiticaT02, KatokN11}, that eliminate
the use of families  and rely only on the \POC.
 When the group is non-commutative, they only require that the cocyles
are close to the identity. This smallness assumption was removed when
$\G$ is finite dimensional in the remarkable paper
\cite{Kalinin11}. As far as we know, results under proximity
to identity (without using families) or, much less, without proximity
assumptions, are not known when the dynamics in the base is
an integrable system as considered here.

\medskip

There are other cases where the cohomology equations
are solvable with estimates (e.g., the base dynamics being a rotation of the torus, interval
exchange transformation, horocycle flow,  twisted versions of the above, higher
rank actions) but these results require other conditions than the 
\POC. In particular, in the case when the base dynamics is a Diophantine rotation,
  the convergence of formal power series expansions in many cases
  was considered in \cite{Moser67}.

Results of the same kind, but with a different condition, exist for partially hyperbolic
  systems, see \cite{Wilkinson13} and references therein.

\bigskip

\subsection{Notations and the Main Theorem}\label{s_notations}

The following notations are needed to give a precise formulation of the results. We will consider  the manifold 
\[
  \M = \torus^d \times [-1,1]^d ,
\]
endowed with the dynamics 
\begin{equation}\label{fdefined} 
  \begin{split} 
    &f: \M \rightarrow \M,  \\
    &f(\th, I) =  (\th + I , I) .
  \end{split}
\end{equation} 

Note that a more general dynamics: $f(\th, I) =  (\th +\Omega( I) , I)$ where  $\Omega$ is invertible, can be reduced to $\Omega(I) = I$ by changing
the variable.  We will refer to $\M$ as the \emph{base} and to $f$ as the
\emph{dynamics on the base}.

\bigskip


Since we are working with the analytic regularity,
it is useful to consider complex extensions of
the above.
Fix $\rho > 0$; we denote:
\begin{equation}\label{notations} 
\begin{split}
& \torus_\rho^d = \{ z \in \complex^d/\integer^d | \, | \text{Im}(z_j)| \le \rho\}, \\ 
  &\B_\rho = \{I \in \complex^d \, | \, \text{dist}(I, [-1,1]^d) \le \rho/10 \},\\
 &\M_\rho = \torus_\rho^d \times \B_\rho, \\
&\E_\rho = \{ \eps \in \complex |\, |\eps| < \rho\}. \\
\end{split}
\end{equation}
We let $\G$ be a Banach algebra (e.g., the algebra of real or complex square matrices or the Banach algebra
of bounded operators in a Banach space).  Later, in Section~\ref{sec:Lie}, we
will present the minimal -- typographical -- modifications needed to adapt the
proofs to $\G$ being a Lie group. 

For any $\rho > 0$, 
we consider $\A_\rho$  to be the space of real-analytic functions on $\M_\rho$, i.e., holomorphic functions
$\phi$ on  
$\M_\rho$, continuous up to the boundary, and satisfying the real symmetry property:
$\overline{\phi (\bar \th,\bar I)}=\phi(\th,I)$ (where the bar stands for the complex conjugate).
In a similar way, we let  $\A_\rho^\eps$ be the space of real-analytic functions on $\cE_\rho \times \M_\rho$.

We endow $\A_\rho$ and $\A_\rho^\eps$ with the supremum norms
 which make  $\A_\rho$ and $\A_\rho^\eps$
into a Banach spaces.  
We keep the same notation, $\| \cdot \|_\rho$, for the two norms.

Since we will need to deal with truncations of power series,
we use the notation: 
 \begin{equation} \label{eq:notation_Taylor}
v^{[m , M]}_\eps = \sum_{j \in [m,M]} \eps^j v^{j} ,
\end{equation} 
 and analogously for other ranges. For example:
 $ 
 v^{[\le M ]}_\eps = \sum_{j \le M} \eps^j v^{j}.
 $

 Notice that $v^{[m , M]}_\eps$  is a polynomial
 of degree $M$ in $\eps$, whose coefficients of order smaller than $m$
 vanish.

 A small typographical confusion is that $\eps^j$ means
 the variable $\eps$ raised to the $j^{\rm th}$ power.
 On the other hand $v^j$ is the $j^{\rm th}$ coefficient
 in the expansion.  This inconsistency  is very common in
 mathematics and we hope will not cause much confusion.

 \bigskip

 We recall the somewhat standard definition of formal power series
 whose coefficients are functions of another variable. Formal
 power series for expansions have been used in mechanics
 and other mathematical disciplines for a long time.

\begin{definition}\label{formalpowerseries}
Given a sequence of continuous                                     
functions $\phi^{j} : \M_\rho \rightarrow \G $,
we say that  expression 
\begin{equation} \label{formal} 
  \phi_\eps = \Id + \sum_{j=1}^\infty \eps^j  \phi^{j}
  \end{equation} 
satisfies
\eqref{coboundary} in the sense of  formal power series
if for any 
$L \in \nat$ we have that the truncated
series $\phi^{[\le L]}_\eps = \Id +  \sum_{j=1}^L \eps^j  \phi^{j} $
satisfies \eqref{coboundary} up to the terms of order $\eps^{L}$:
\begin{equation}\label{truncatedbounds}
\| \eta_\eps  - (\phi_\eps^{-1})^{[\le L]} \circ f\cdot  \phi_\eps^{[\le L]}\|_{C^0} \le C_L |\eps|^{L+1}.
\end{equation}

\end{definition}

We emphasise that \eqref{formal} is just a suggestive formal expression 
and that it is not assumed to converge in any sense.  The only meaning
adscribed to the infinite sum is precisely that each of the truncated sums satisfies
\eqref{truncatedbounds}.

Of course, we could consider other spaces rather than $C^0$
for the functions
in the coefficients, or even spaces depending on the order
(a common situation in series expansions in mechanics
is that the coefficients are analytic, but the analyticity domain
depends on $n$ and even shrinks as $n$ tends to $\infty$). 

\bigskip

Here is the main result of the paper.
\begin{theorem}\label{main}
  Assume the notation in \eqref{fdefined}, \eqref{notations}, and consider an analytic family 
  $\eta_\eps:  \M_\rho \rightarrow \G$, $\eps \in \cE_{\rho}$,
 with $\eta_0 = \Id$. 

  Then the following statements are equivalent:
  \begin{itemize} 
  \item[$(A)$.]  $\eta_\eps$ satisfies the \POC\  \eqref{POC-general} for each $\eps\in \cE_\rho$.
  \item[$(B)$.]There is a solution of
    \eqref{coboundary} in the sense of formal power series expansions
    (see Definition~\ref{formalpowerseries}).
  \item[$(C)$.] There exists a real-analytic solution $\phi_\eps$ of
    \eqref{coboundary}. The domain of analyticity of $\phi_\eps$ may be
    smaller than the domain of $\eta_\eps$.
  \end{itemize}
\end{theorem}

\bigskip

\subsection{Comments on the proof}\label{s_comments_on_pf}

Clearly, $(C) \implies (A), (B)$.
The fact that $(B) \implies (A)$ is also easy.
We observe that for any periodic point $p$ of $f$,
the function $\eta_\eps(f^{N-1}p) \cdots \eta_\eps(p)$ is
analytic in $\eps$. The existence of a formal power series solution of
\eqref{coboundary} implies that all the terms in the power
series expansion vanish.

The fact that $(A)$ implies  $(C)$ is the main result of this paper. It will be proved in
 Section~\ref{sec:proof} using a Nash-Moser iterative approach.
  The basic tool in this proof is solving the  linearized equation with estimates. 
Namely, by linearizing the coboundary 
equation (using the fact that $\eta_\eps$ is close to the identity), we obtain a 
commutative cohomology  equation over 
$f$. We show that if the commutative cohomology equation has a formal power series solution,
then it has an analytic solution  with good estimates 
(the estimates involve a loss of the domain with
respect to the initial data). 
Note that
due to the lack of uniqueness, the solution with estimates may be
different from the formal power solution assumed in the hypotesis.

Using the group structure of the coboundary equation and
of the hypothesis, we can use this approximate solution
to improve the situation.  The process can be repeated
and the convergence of this accelerated procedure is
established by a slight  modification of  the usual Nash-Moser
scheme (we get an extra term due to the truncation of
power series).


\bigskip

\subsection{Various remarks}

\begin{remark} 
Notice that we are not stating 
that the formal power series solution in the hypothesis $(B)$ of
Theorem~\ref{main} 
converges (indeed, such result is false!).  

In fact, the solutions of
the coboundary equation are far from being unique.  Indeed,  if $\phi_\eps(\th,I)$ satisfies coboundary equation \eqref{coboundary} 
and $A_\eps (I)$ is a function from  $\E_\rho \times B_\rho$ to $\G$,
then
\begin{equation} \label{nonunique}
  \tilde \phi_\eps(\th, I) =  A_\eps (I)\cdot  \phi_\eps(\th, I)
 \end{equation} is also
a solution of the coboundary equation.

We show
that, from the existence of  a formal power series solution
we can obtain a (possibly different) power series
that converges in some complex domain. 
For this proof it is important that we work with
analytic one-parameter families of cocycles.
 \end{remark}

\begin{remark} 
  For simplicity of notation, we have used only one parameter, $\rho$,
  to indicate the size of the analyticity domains in $\eps, I, \th$.
  Since  these variables play different roles, it would have
  been natural, and improve slightly the estimates, to consider them as independent. 
  For present purposes the improvement in
  regularity does  not justify the complication of dealing with norms indexed by 3 parameters. 
  \end{remark} 

\begin{remark}
  The reason to take $\G$ to be a Banach algebra is
  that it allows to apply corrections 
by adding terms. This simplifies the notation, but is not
essential. One can consider $\G$ to be a Lie group, the correction
functions taking values on the Lie algebra, in which case the corrections are
applied using the exponential mapping rather than just adding. After
the proof  of Theorem~\ref{main} is presented, we collect in Section~\ref{sec:Lie} the (mostly
typographical) changes needed to obtain the result for $\G$ being a
(Banach) Lie Group.
\end{remark} 

\begin{remark} \label{realvalues} 
Note that $f$ does not map the  complex domain $\M_\rho$ into itself and this will be a
source of -- mainly notational -- problems. We have
that $f(\M_\rho) \subset \M_{11\rho/10}$ (when $I$ is complex, an application of $f$
changes the imaginary part of the $\th$ component.  This
is the reason for introducing the factor $10$ in the definition of
$\B_\rho$).
\end{remark}

\begin{remark} 
Note that \eqref{POC-general} is far from being a trivial condition.
The map $f$ indeed has a lot of periodic orbits. If there exists $N\in \nat$ such
that $N I \in \integer^d$, then for any $\th \in \torus^d$, $(\th, I)$ is periodic of
period $N$.  Note also that the periodic orbits of $f$ are real, hence POC 
only provides information about the behaviour of $\eta_\eps$ in restriction to the real hyperplane.
It appears that this information is sufficient in order to obtain an estimate of the solution in a complex neighborhood. 
We would call this a manifestation of the ``magic" of complex analysis.
 \end{remark}

 \begin{remark} \label{nolower}
    One
  can wonder if it would be possible to consider
  analogous  results when  the dimension of
  the angles $\theta$ in the base is bigger than the dimension of
  the actions $I$. The simple example
  $f(\theta_1,\theta_2, I) =  (\theta_1 + \omega, \theta_2 + I) $
  shows that such dynamics in the base can lack periodic orbits
  so that the \POC\
  becomes vacuous (hence, trivially satisfied) 
  but not all cocycles are coboundaries.
One can ask whether, in this context,
  the existence of 
  formal power solutions implies a convergent solution. Instead of the periodic orbit condition one can explore 
  more general obstructions: invariant
  measures, pseudomeasures as in \cite{Veech86} or distributions.
  \end{remark}


  \section{Commutative cohomology equations over integrable  dynamics}

The key tool we will use in the proof of Theorem~\ref{main}
consists in estimates of the solutions to a commutative version of
the coboundary equation in the case  when
the group $\G$ is commutative.
Namely, in this section we will consider equations for $\alpha: \M_{\eps} \rightarrow \real^n$
given $\beta: \M_{\eps}  \rightarrow \real^n $
of the form:
\begin{equation} \label{commutative_equation}
  \alpha \circ f - \alpha = \beta .
\end{equation}

Notice that the solution to \eqref{commutative_equation}
is not unique. For instance, any function $\alpha$ only depending on $I$ (and not on $\th$) 
leads to $\beta = 0$.
This lack of uniqueness of solutions to the commutative equation can be seen as a
commutative version of the non-uniqueness \eqref{nonunique} for
the general equation. 

The equation will be solved using Fourier series in the angle variable of the following form:
\begin{equation}\label{partialFourier}
\beta(\th, I) = \sum_{k \in \integer^d} \hbeta_k(I) e^{2 \pi i \<k,\th \>} 
\end{equation}
(and analogously for other functions).
With  this notation, under very small regularity
requirements on $\alpha$, we have 
\[
  \alpha \circ f(\th, I) = \sum_{k \in \integer^d} \halpha_k(I)
  e^{2 \pi i \<k,I \> }  \, \cdot \, 
  e^{2 \pi i \<k,\th \>}, 
\]
and \eqref{commutative_equation} is equivalent to the
following set of equations (for $\{ \halpha_k(I) \}_{k \in \integer^d} $ given 
 $\{ \hbeta_k(I) \}_{k \in \integer^d}$):
\begin{equation}\label{commutative_equation_Fourier}
  \halpha_k(I) \left( e^{2 \pi i \<k,I \>} - 1  \right)
  = \hbeta_k(I).
  \end{equation}

The theory we will develop will consist in:
\begin{itemize}
\item Characterization of obstructions for solvability; 
\item Characterization of lack of uniqueness of solutions; 
\item Estimates for the solutions (when they exist).
\end{itemize}

\subsection{Obstructions CPOC and FC for the existence of
  solutions of \eqref{commutative_equation}}

Here we present two necessary conditions for the existence
of solutions of \eqref{coboundary}.
The first one is  based on the study of the dynamics, and
the second one uses Fourier representations.
We will show that they are actually equivalent.
As we will see, the geometric representation is useful when we
perform changes of variables.  The Fourier representation
is useful to obtain estimates of the solutions.
Combination of these two points of view  will allow us to make progress. 

\subsubsection{Periodic Orbit Condition for the commutative equation (CPOC)}

Suppose that  equation \eqref{commutative_equation} has a solution.
If $p = f^N(p)$ is a periodic point for the base dynamics, 
 then the sum of the values of  $\beta$ over
the periodic orbit telescopes to zero:
\[
  \sum_{j= 0}^{N-1}  \beta\circ f^j(p) = 0.
\]
In our case for each $N\in\NN$ we have $f^N(\th, I_*) = (\th + N I_*, I_*)$.
Hence, $(\th, I_*)$ is periodic of period $N$ if and only if $N I_* \in \integer^d$. 
This motivates the following definition.

\begin{definition} We say that  a function $\beta (\theta ,I)$ satisfies the Commutative Periodic Orbit Condition (CPOC)
if  for every $I_*\in \B_\rho$ satisfying
 $N I_* \in \integer^d$   for some $N \in \nat$ we have:
\begin{equation}\label{eq:POC}
\sum_{j=0}^{N-1} \beta(\th + j I_*, I_*)  = 0 \quad \text{for all real } \theta.
\end{equation}
\end{definition}

\begin{remark} \label{rem_CPOC_real}
Note that the $I_*$ that give rise to periodic orbits are
of the form $\ell/N$ with $\ell \in \integer^d$. In particular,
such $I_*$ are real. This will lead to some difficulties. 
\end{remark}

The following result summarises the discussion in the beginning of the section.

\begin{proposition}\label{prop:POC} 
Let  $\alpha$ and $ \beta$ be continuous
functions satisfying \eqref{commutative_equation}. 
Then $\beta$ satisfies CPOC.
\end{proposition}

\subsubsection{Fourier coefficient condition (FC)}

\begin{definition} We say that  a function $\beta (\theta ,I)$ satisfies Fourier Coefficient condition (FC)
if for every  $k \in \integer^d \setminus \{0\}$ and every $I_*\in \B_\rho$ we have:
  \begin{equation}\label{eq:FC}
\< k,  I_* \> \in \integer \ \text{ implies }  \ \hbeta_k(I_*) = 0.
  \end{equation}

\end{definition}

Looking at \eqref{commutative_equation_Fourier} we
obtain immediately the following result.

\begin{proposition}
  A necessary condition for the existence of
  continuous solutions of \eqref{commutative_equation} is
  that  $\beta$ satisfies FC.
\end{proposition}

The above statement tells us that, in order  for the commutative equation \eqref{commutative_equation}
to have a continuous solution $\alpha$,
the Fourier coefficient
of every index $k$ has to vanish for $I_*$ lying in parallel planes
corresponding to $k$. Note, however, that these planes involve complex values. 
In contrast, condition CPOC gives us information only about real values of $I$, see 
Remark \ref{rem_CPOC_real}.
This makes the following result rather surprising.

\subsubsection{Equivalence of conditions FC and CPOC}

Conditions FC and CPOC
have  a very different nature. CPOC is
geometrically natural, and it is clearly preserved under
a cocycle conjugacy. FC, on the other hand, leads
to very effective estimates of
the solution, as we will see in Section~\ref{sec:estimates-commutative}. Below we prove that the two conditions  are equivalent.

\begin{proposition} \label{POC-FCequivalence} 
A continuous map $\beta: \M_\rho \rightarrow \real$ satisfies  CPOC (condition 
  \eqref{eq:POC}) if and only if it satisfies FC (condition \eqref{eq:FC}).
\end{proposition}

 At the first glance, CPOC may seem weaker than FC, since the it
gives us 
information only about the real values of the argument $I$,
while  FC concerns a complex neighbourhood in $B_\rho$. 
The following uniqueness  theorem from complex analysis will
be used to bridge this gap.

\begin{lemma}\label{wonderful_argument}
If $f(I)$ is an analytic function on a domain $D\in \CC^d$ that vanishes in a real neighbourhood $U^{\RR}$ of a point $I_0=x_0+iy_0\in D$, that is, on a set 
$$U^{\RR}=\{I=x+iy\in \CC^d \mid \, |x-x_0|<r,\ y=y_0\},$$ 
then $f(I)\equiv 0$ on $D$.
\end{lemma}
The idea of the proof is very simple. We note that all the (complex) derivatives
of the function can be computed along the real space. The assumption implies
that the real derivatives vanish. Hence, all the derivatives vanish. 

\begin{proof}
Since $f(I)\equiv 0$ on $U^{\RR}$, for any  $v\in \RR^d$, for $t\in\RR$ and $j\in\NN$ we have 
$\left(\frac{d}{dt}\right)^j f (u+tv)\mid_{t=0}=0$. Hence, 
$$
\left(\frac{d}{dt}\right)^j f (u+tv)\mid_{t=0}=
\left( \<v, \nabla \> \right)^j f (u)=0.
$$
Since this is true for all  $v\in \RR^d$, we get $\partial^k f (u)=0$ for all multi-indices $k\in \NN^d$. 
Since $f$ is analytic in $D$, this implies that $f$ vanishes identically.
\end{proof}

\begin{proof}[Proof of Proposition \ref{POC-FCequivalence}.]
The direction of the equivalence that is used in this paper
is that CPOC  
implies FC, 
so we do it first.

Suppose that CPOC (condition \eqref{eq:POC}) holds. 
Fix any $k \in \integer^d$, $n\in \ZZ$, and consider the complex hyperplane 
$P^{\CC}_{k,n}=\{I \in \CC^d\mid \<k,I \> = n \}$. Denote $P^{\RR}_{k,n} =\{I \in P^{\CC}_{k,n} \mid  \Im I=0 \}$.

Note that, since $k\in\ZZ^d$, the points of
the form $\ell/N$ with $\ell \in \integer^d$ and $N\in \nat $ are dense on  
$P^{\RR}_{k,n}$. 
Fix any $I = \ell/N\in P^{\RR}_{k,n}$. 
 Using that $\<  k, t \ell \> \in \integer$, we get: 
  \[
    \begin{split} 
    \hbeta_k(\ell/N) & =
    \int_{\torus^d} \beta(\th, \ell/N) e^{-2 \pi i \< \th , k \> } \, d\th  \\
    &  = \int_{\torus^d} \beta(\th + \ell/N , \ell/N) e^{-2 \pi i \< \th , k \>  }
   e^{-2 \pi \< \ell/N , k \>} \, d\th  \\
  & = \int_{\torus^d} \beta(\th + \ell/N , \ell/N) e^{-2 \pi i \< \th , k \>  }  \, d\th\\
   &= \cdots \\
   &=\int_{\torus^d} \beta(\th + j  \ell/N , \ell/N) e^{-2 \pi i \< \th , k \>  }  \, d\th  . \\
    \end{split}
  \]
 Adding the above $N$ different expressions for $\hbeta_k (\ell/N)$, we obtain:  
  \begin{equation}\label{identity_on_periodic} 
    N \hbeta_k (\ell/N)  =   \int_{\torus^d}
    \left( \sum_{j = 0}^{N-1}
      \beta(\th + j  \ell/N , \ell/N) \right)  e^{-2 \pi i \< \th , k \>  }  \, d\th  . \\
  \end{equation} 
Under the CPOC condition, the integrand in \eqref{identity_on_periodic}
is identically zero. We have thus shown that $\hbeta_k(\ell/N)=0 $ for any rational point $(\ell/N)\in P_{k,n}^\RR$.
Since such points are dense on $ P_{k,n}^\RR$,  
and the zero-set of a continuous function is always closed,   we have:
$\hbeta_k(I) =0$ for any  $I\in P_{k,n}^\RR$.  Since the restriction of $\hbeta_k(I)$ onto $P^{\CC}_{k,n}$ is analytic, 
the above result, combined with Lemma \ref{wonderful_argument}, implies that  the 
restriction of $\hbeta_k(I)$ onto $P^{\CC}_{k,n}$ equals zero.
Hence, FC condition \eqref{eq:FC} holds.

  \bigskip

  Now we prove that FC implies CPOC. Assume that FC holds. Let $I=\ell/N $ for some $\ell \in \integer^d$ and $N\in \nat $. Then
\eqref{identity_on_periodic} holds.  We note that the function
  \[
  A(\th)   = 
\left( \sum_{j = 0}^{N-1}
\beta(\th + j  \ell/N , \ell/N) \right)
\]
satisfies $A(\th + \ell/N) = A(\th)$, so it is a function defined
in a reduced torus obtained by identifying the points that differ by a translation by $\ell/N$.
For functions with this extra  periodicity, the exponentials $e^{2\pi i \<k,\theta \>}$
with  condition $\<k,\ell /N \> \in \ZZ$ form a complete set.
Therefore, from  FC,  we conclude that $A(\th) \equiv 0$ and, hence,
that  the CPOC condition \eqref{eq:POC}   is satisfied.
\end{proof}

\begin{remark} 
The fact that
FC implies CPOC is not used in the proofs below. We can get this statement also as a
 byproduct of
the main line of argument. 
Indeed, in the next section we will show
that for analytic  functions $\beta$ satisfying
FC, one can construct an (analytic) solution $\alpha$ to
\eqref{commutative_equation}. The fact that equation \eqref{commutative_equation} 
with the right-hand side $\beta$ has a continuous solution
implies  that $\beta$ satisfies CPOC.
\end{remark}

\subsection{Estimates for  the commutative equation \eqref{commutative_equation}}\label{sec:estimates-commutative}
\subsubsection{A reminder on Cauchy estimates} Assume the notations from  Section  \ref{s_notations}; in particular consider the space $\cA_\rho$ of analytic functions.
The following statements are standard. 
   \begin{lemma}\label{lem:Cauchy}
     If $v\in \cA_\rho$,  $v(\th,I)=\sum_{k\in \ZZ^d} v_k (I) e^{2\pi i \< \th , k \> }$,
    then there exists a constant $c=c(d)$, such that for each $I\in \B_\rho$ we have:
\begin{equation} \label{eq:Cauchy}
  \begin{split}
   & |v_k|\leq \|v\|_{\rho} e^{-2\pi |k| \rho }, \\
    &\| D v \|_{\rho - \delta} \le c \delta^{-1} \| v \|_{\rho}. \\
     \end{split}
\end{equation} 
Given  $v_\eps \in \cA_\rho^\eps$, write it as a Taylor series in $\eps$: $\displaystyle v_\eps(\th, I)=\sum_{j=0}^{\infty}\eps^j v^{j}(\th, I)$.
Suppose that $\rho$ satisfies $0< \rho\leq\rho_0$ for some fixed $\rho_0$.
Then, recalling notation~\eqref{eq:notation_Taylor}, we can estimate: 
    \begin{equation} \label{eq:Cauchy2}
  \begin{split}
  & \| v^{j} \|_{\rho}\leq  \rho^{-j}  \|v\|_{\rho} , \\
    & \| v^{[m , M]} \|_{\rho -\delta} \leq \rho_0 \delta^{-1}  \|v\|_{\rho} , \\
     & \| v^{[\geq M]}\|_{\rho -\delta} \leq \rho_0  \delta^{-1}  \exp ({- M \delta/\rho_0 }) \|v\|_{\rho_0} . \\
  \end{split}
\end{equation}

 \end{lemma}  
 
 \begin{proof}  Let us present the proof of the last estimate for completeness. First of all, since $\ln (1-x)\leq -x$ for all $x\in \RR$, we have
$$
(1-\dt/\rho)  \leq  e^{-\dt/\rho}.
$$ 
 Using $ \| v^{j} \|_{\rho}\leq  \rho^{-j}  \|v\|_{\rho}$ and $\eps <\rho-\dt$, we have:
\[  
\begin{aligned}
& \| v^{[\geq M]}\|_{\rho -\delta} \leq   \sum_{j=M }^{\infty} \eps^j \| v^{[j]} \|_{\rho} \leq  \sum_{j=M }^{\infty} (\rho-\dt)^j \rho^{-j}  \|v\|_{\rho} \\
 &\leq   \|v\|_{\rho} (1-\dt/\rho)^M \sum_{j=0 }^{\infty} (1-\dt/\rho)^j  \leq  
\|v\|_{\rho} \,    e^{- M \delta/\rho_0  } \,      \rho \delta^{-1} \\ 
& \leq \rho_0  \delta^{-1}  e^{-  M \delta /\rho_0 } \|v\|_{\rho} . \\
  \end{aligned}
\]   
 
\end{proof}
      \subsubsection{Estimates on the solutions of the  coboundary equation \eqref{commutative_equation} }

      \begin{lemma}\label{commutative_estimates}
        Assume that $\beta $ is analytic in $\M_\rho=\torus^d_{\rho}\times \B_\rho$ and satisfies the
        FC condition \eqref{eq:FC}.

        Then there exists $\alpha$ solving 
        equation~ \eqref{commutative_equation}, i.e.,
$\alpha \circ f - \alpha = \beta$,
 such that, for a constant $c=c(d)$ and for any $0 < \delta < \rho$ we have:
        \begin{equation} \label{eq:tame}
          \|\alpha \|_{\rho - \delta} \le c \delta^{-d-1} \|\beta\|_\rho.
        \end{equation}

        If  $\beta$ depends analytically (continuously)  on a
        parameter $\eps$ ranging in a certain domain, a solution
        $\alpha_\eps$ can be chosen to depend analytically (continuously) on
        $\eps$ in the same domain. 
      \end{lemma}

   \begin{remark}   As an important corollary of this lemma, we obtain estimates of the solution in domains
      other than $\M_\rho$.
      Namely, using equation~\eqref{commutative_equation}, we obtain (for a slightly modified constant, which we denote by $c$ again):
        \begin{equation} \label{eq:tame2}
          \|\alpha\circ f  \|_{\rho - \delta} \le c \delta^{-d-1} \|\beta\|_\rho.
        \end{equation}
\end{remark}

      \begin{proof}
       Equation \eqref{commutative_equation} is equivalent to
        the sequence of equations \eqref{commutative_equation_Fourier}
        for the partial Fourier coefficients
        defined in \eqref{partialFourier}:
  \begin{equation}\label{commutative_equation_Fourier2}
  \halpha_k(I)\left( e^{2 \pi i \< k,  I\>} - 1  \right)
  = \hbeta_k(I).
\end{equation}
For $I$ such that $\< k,  I\> \notin \integer$ we can express $  \halpha_k(I)=\hbeta_k(I) \left( e^{2 \pi i \< k,  I\>} - 1  \right)^{-1}$.
The crucial remark is that the FC condition for $\beta$ implies that 
$$
\< k,  I\> \in \integer \Rightarrow \hbeta_k(I)=0.
$$
Hence, for  $I$ such that  $\< k,  I\> \in \integer$,
equation \eqref{commutative_equation_Fourier2} is satisfied for any value of
$\halpha_k(I)$. 
We define $\halpha_k(I)$ for these $I$ by continuity. 
A way to do it is the following (compare with the L'H\^opital's rule). 
Differentiate equation \eqref{commutative_equation_Fourier2} 
in the direction of the vector $k$:
$$
\< \nabla_I  \halpha_k(I) , k \> \left( e^{2 \pi i \< k,  I\>} - 1  \right) + 
\halpha_k(I)   2 \pi i \< k,  k\> e^{2 \pi i \< k,  I\>}
= \<\nabla_I \hbeta_k(I), k \> .
$$
If $\< k,  I\> \in \integer$, this gives us $2 \pi i | k|^2 \halpha_k(I)    e^{2 \pi i \< k,  I\> }
= \<  \nabla_I \beta_k(I) , k \>$.
Summing up, we have defined a continuous function $\halpha_k(I)$ by     
\begin{equation}\label{commutative_equation_Fourier3}
  \halpha_k(I)
  = \begin{cases}
     \hbeta_k(I)\left( e^{2 \pi i \< k,  I\>} - 1  \right)^{-1}
   & \text{if   } \< k,  I\> \notin \integer \\
  \<  \nabla_I \beta_k(I) , k \>  e^{-2 \pi i \< k,  I\> } / ( 2 \pi i | k|^2 )& \text{if   } \< k,  I\> \in \integer . \\
    \end{cases} 
\end{equation}
Since $\halpha_k(I)$ is analytic in $\B_\rho \setminus \{\< k,  I\> =0\}$ and
bounded in $\B_\rho$, it is analytic in $\B_{\rho}$.

Now let us estimate the norm of the solution. This estimate is done in \cite{dlLS}, Lemma 6; below we repeat the argument
for completeness. Fix $0<\dt<\rho/2$.
  For each fixed 
$k\in\ZZ^d$,  we will estimate the corresponding $\halpha_k(I)$ in two steps: first `$\delta/2$-close" to the resonant plane $\< k,  I\> $, and then in the rest of $\B_{\rho -\dt} $. 

For the first step, let $\Pi_\dt= \{\< k,I\>=0 \} \cap \B_{\rho-\dt}$ be the part of the 
resonant plane falling into $\B_{\rho-\dt}$. Notice that the orthogonal complement to this plane is formed 
by the vectors $\gamma e^{2 \pi i \theta } k $, $\gamma \geq 0$, $\theta \in [0,1)$. Let 
$$
\Delta=\left\{ I= \gamma  \frac{ k }{| k | } e^{2 \pi i \theta } \  \Big| \  \gamma < \dt  /2, \,\theta \in [0,1) \right\}
$$ 
be the complex disk 
of radius $\dt /2$ centered at zero and orthogonal to $\Pi_\dt$.  
Note that the restrictions of
$\halpha_k(I)$ and $\hbeta_k(I)$ to this disc 
are analytic.
Consider the ${\dt}/2$-neighbourhood $O_\dt$ of 
$\Pi_\dt$: 
$O_\dt=\bigcup_{I_0\in \Pi_\dt} (I_0+\Delta)$. Then $O_\dt\subset \B_{\rho -\dt} $.

For each fixed $I\in O_\dt$ there exists $I_0 \in \Pi_\dt$ such that $I\in I_0+\Delta$.  We can estimate $|\alpha_k (I)|$ by
the maximum modulus principle on the disk  $I_0+\Delta$. Namely, for $I$ lying on the boundary of this disk we have: $|\< k,  I\>|= 
|\< k,  I_0\>+ \<  k, \dt  k/(2| k|) \> |  = | k|\dt/2$. Hence, for such $I$ we have
$$
|\halpha_k (I)| \leq 
\frac{2 \| \hbeta_k \|_\rho }{4\pi\dt | k|}< \frac{ \|\hbeta_k \|_\rho }{\dt | k|}.
$$
As the second step in this estimate, consider $I\in \B_{\rho -\dt}\setminus O_\dt $. Here 
$|\< k,  I\>| \geq  |k| \dt / 2$, so $|\halpha_k (I)| $ satisfies the same estimate as above. This proves the desired estimate for an individual parameter value.

In the case that the data depend analytically (continuously)
on a parameter $\eps$, we obtain that the Fourier
coefficients depend analytically (continuously) on 
$\eps$. The uniform bounds on the Fourier coefficients imply
that the sum depends analytically (continuously) on $\eps$. 
    
\end{proof} 
        
      \section{Proof of Theorem~\ref{main}}
      \label{sec:proof}

Theorem ~\ref{main} is the equivalence between  the existence of an analytic solution to the coboundary equation, $(C)$, and  two formal conditions:
\POC\  $(A)$,
and the existence of formal power series solutions, $(B)$.
As indicated in Section \ref{s_comments_on_pf}, the heart of
the problem is to show that $(A)$ implies $(C)$. 

\subsection{Overview of the proof}
\label{sec:overview}
The main part of the proof consists in estimating  the results of
an iterative step.

For each $n\in \ZZ$, we let $L_n=L_0 2^n$, where $L_0$ is an appropriate constant. 
The initial input of the iterative step  will be an
$\eta^n_\eps = \Id + O(\eps^{L_n})$
satisfying the POC condition \eqref{POC-general}.
The iterative procedure will produce an ``almost solution" $\phi^{n}_\eps = \Id + O(\eps^{L_{n}})$
such that
\begin{equation}\label{improved} 
  \eta^{n+1}_\eps  :=  \left(\phi^n_\eps\circ f \right)^{-1} \cdot \eta^n_\eps  \cdot \phi^n_\eps
\end{equation}
satisfies $\eta^{n+1}_\eps = \Id + O(\eps^{2L_{n}})$. Note that, because
of the construction of \eqref{improved}, we obtain that
$\eta^{n+1}_\eps$ also satisfies the POC condition,
and the iterative procedure can be applied again.

We will show that the iterative procedure leads to the estimates
\[
\|\eta^{n+1}_\eps - \Id\|_{\rho_n - \delta_n}
\le c \delta^{-2(d+1)} \|\eta^{n}_\eps - \Id\|_{\rho_n}^2 +  c\delta_n^{-1}\|\eta^{n}_\eps - \Id\|_{\rho_n} 
e^{- L_n \frac{\delta_n}{\rho_0}},
\]
as well as estimates for $\phi^n_\eps$.
From these estimates, using standard arguments in Nash-Moser theory (presented in Section~\ref{sec:convergence}),
we will conclude
that the limit 
$\displaystyle \lim_{n \to \infty} \phi_\eps^1\cdot \phi_\eps^2  \cdots \phi^n = \phi_\eps^\infty
$
exists in a domain, and that
\[
(\phi^\infty_\eps\circ f)^{-1} \cdot \eta_\eps\cdot \phi^\infty_\eps = \Id .
\]

\subsection{Formal construction of the iterative step}\label{s_iter_step_formal}
In this section, we will describe the formal procedure leading to an improved solution; the estimates are provided in the next section. 
Given an analytic function $\eta_\eps(\theta,I)$ and $L\in\mathbb N$, we use  notation \eqref{eq:notation_Taylor}
to write
 \begin{equation}\label{highorderform}
  \eta_\eps = \Id + \eta_\eps^{[L, 2L -1]} + O(\eps^{2L} ),
 \end{equation}
 so that $\eta_\eps^{[L, 2L -1]}$ is a polynomial of
 degre $2L -1$ in $\eps$, but its coefficients of order smaller
 than $L$ vanish. 

\begin{proposition}\label{POC->CPOC} 
Suppose that $\eta_\eps$ has the form \eqref{highorderform} and satisfies the POC condition
\eqref{POC-general}.
Then the truncated function,
$\eta_\eps^{[L, 2L -1]}$, satisfies the CPOC condition
\eqref{eq:POC}.
\end{proposition}

\noindent{\it Proof. \ } Given $\eta_\eps$ as above 
and a periodic orbit $p = f^N(p)$,  note that
\begin{equation}\label{expansion_periodic} 
\eta_\eps( f^{N -1}(p) ) \cdots \eta_\eps(p) =
\Id +  \sum_{j = 0}^{N-1} 
\eta_\eps^{[L, 2L -1]}( f^j(p)) + O(\eps^{2L}).
\end{equation}
By POC condition, the left hand side equals $\Id$ for all $\eps$. Note that $\eta_\eps^{[L, 2L -1]}$ is a polynomial in $\eps$ (of degree $2L-1$). Therefore, the sum in the right hand side equals zero, which is precisely the CPOC condition
for $\eta_\eps^{[L, 2L -1]}$. \qed

\begin{remark}
Notice that Proposition~\ref{POC->CPOC}  depends crucially on the fact that we
are using families of maps, and assume the Periodic Orbit Condition for all values in $\eps$. In fact, this is the main reason for the use of families in this paper. Note also that the approximation given in \eqref{expansion_periodic} is very non-uniform in $N$. 
\end{remark}

By Proposition \ref{POC-FCequivalence}, condition CPOC for $\eta_\eps^{[L, 2L -1]}$ implies condition FC for $\eta_\eps^{[L, 2L -1]}$.  

Under the FC condition for $\eta_\eps^{[L, 2L -1]}$, Lemma~\ref{commutative_estimates} gives us 
a  $\phi^{[L, 2L -1]}$ solving
\begin{equation}\label{iterativesol}
 \phi^{[L, 2L -1]}_\eps \circ f 
- \phi^{[L, 2L -1]}_\eps = \eta^{[L, 2L -1]}_\eps .
\end{equation}

It is easy to see that, defining $\phi_\eps = \Id + \phi_\eps^{[L, 2L -1]}$,
  we have $\phi_\eps^{-1} = \Id - \phi_\eps^{[L, 2L -1]} +O(\eps^{2L})$. For 
$\eta_\eps=\Id +\eta_\eps^{[\geq L]}$ we have:
 \begin{equation}\label{eta-id-formal}
\begin{aligned}
  \tilde \eta_\eps - \Id: =&\phi_\eps^{-1}\circ f  \cdot  \eta_\eps  \cdot \phi_\eps  -\Id = \\
 =&   ( \Id + \phi^{[L, 2L -1]}_\eps \circ f)^{-1} \ \cdot (\Id +\eta_\eps^{[\geq L]})  \cdot ( \Id + \phi^{[L, 2L -1]}_\eps ) -\Id \\
=& A+ (  \Id - \phi^{[L, 2L -1]}_\eps \circ f) \cdot (\Id +\eta_\eps^{[\geq L]}) \cdot  ( \Id + \phi^{[L, 2L -1]}_\eps )  -\Id  \\
 =& A+  B + \eta_\eps^{[\geq L]}  + \phi^{[L, 2L -1]}_\eps - \phi^{[L, 2L -1]}_\eps \circ f  \\
=& A+  B + C + \eta^{[L, 2L -1]}_\eps  + \phi^{[L, 2L -1]}_\eps - \phi^{[L, 2L -1]}_\eps \circ f  \\
=& A+  B + C ,
\end{aligned}
\end{equation}
where
$$
\begin{aligned}
&\tilde A= ( \Id + \phi^{[L, 2L -1]}_\eps \circ f)^{-1} -   (\Id - \phi^{[L, 2L -1]}_\eps \circ f),  \\
&A= \tilde A \cdot \eta_\eps  \cdot ( \Id + \phi^{[L, 2L -1]}_\eps )  ,  \\
&B=    \eta_\eps^{[\geq L]}  \cdot \phi^{[L, 2L -1]}_\eps - \phi^{[L, 2L -1]}_\eps \circ f  \cdot (\eta_\eps^{[\geq L]} +
\phi^{[L, 2L -1]}_\eps +
\eta_\eps^{[\geq L]} \cdot \phi^{[L, 2L -1]}_\eps), \\
&C = \eta^{[\geq  2L] }_\eps,
\end{aligned}
$$
and the last equality in \eqref{eta-id-formal}  holds since $\phi^{[L, 2L -1]}_\eps$ is the solution of equation \eqref{iterativesol}.
\qed

 \medskip
  

 \begin{remark} \label{r_eta_tilde}    Note that expression above readily implies that 
$$
\tilde \eta_\eps=\phi_\eps^{-1}\circ f  \cdot \eta_\eps  \cdot  \phi_\eps=\Id + O(\eps^{2L}).
$$
Hence, the formal argument developed here shows  that the
      POC condition (i.e., condition $(A)$ in Theorem ~\ref{main})
      implies the existence of a formal power series solution
      $(B)$ in Theorem ~\ref{main}). 
      \end{remark} 
To prove that the resulting limit function $\phi_\eps^{\infty}$ is a convergent power series, 
we need to develop detailed
estimates that establish the convergence of the iterative procedure. The key step is
to show that the error in the procedure decreases faster than exponentially. This will show
that the corrections applied to $\phi$ are summable, so that there is a limit. 

The estimates are somewhat delicate because step
\eqref{iterativesol} involves a loss of domain with singular
estimates.  So, to control the convergence we will use the arguments from  KAM theory.

\subsection{Estimates on the iterative step}
In this section we repeat the iterative step described in the previous section, 
while keeping track of the sizes of all the objects involved.  The main technical tool
is Lemma~\ref{commutative_estimates}. We prove the following.

\begin{proposition} \label{prop:iterative} Given numbers $0 < \rho\leq \rho_0$, $L\in \NN$ and $\sigma>0$,  
let $\eta_\eps\in \cA_\rho^\eps $ be  an analytic family taking values
  in $\G$, such that  $\eta_0 = \Id$, and for each $\eps\in \cE_{\rho}$ we have:
  \begin{itemize}
  \item
$
\eta_\eps =
  \Id + \eta^{[L, 2L -1]}_\eps + \eta_\eps^{[\ge 2L]},
$

 \item
    $\eta_\eps$ satisfies the POC condition \eqref{POC-general},
\item $ \|\eta_\eps - \Id\|_\rho  \le \sigma$.
  
 \end{itemize}

  Then we can find $\phi_\eps = \Id + \phi_\eps^{[L, 2L -1]}$, 
  for  $s=d+1$ and a certain constant $c=c(d)$  satisfying estimates
 $$ \|\phi_\eps^{[L, 2L -1]}\|_{\rho -\delta}    \le c \delta^{-s} \sigma 
 ,\quad
    \|\phi_\eps^{[L, 2L -1]} \circ f \|_{\rho -\delta}    \le  c \delta^{-s} \sigma,
$$
 such that the following holds. 
  Define
 $$  
  \tilde \eta_\eps := \phi^{-1}_\eps \circ f \cdot \eta_\eps\cdot  \phi_\eps .
$$
Then
  \begin{itemize}
 
 \item
    $ \tilde \eta_\eps = \Id +  \tilde\eta_\eps^{[\ge 2L]}$,
  
  \item
    $\tilde \eta_\eps$ satisfies the POC condition
    \eqref{POC-general}.
   
  \item
Assuming furthermore that, for a certain constant $c'=c'(d)$,  $\delta$ satisfies
  \begin{equation}\label{iterativeassumption}
  c'  \delta^{-s}\sigma  \le 1/100,
  \end{equation}
  we have that $\tilde \eta_\eps $   satisfies:
 \begin{equation}  \label{est-C}
  \| \tilde \eta_\eps - \Id \|_{\rho - \delta}
    \le  c \delta^{-2 s} \sigma^2 + c\delta^{-1}  \si \exp{( -  \delta L/\rho_0)} .
\end{equation} 
 
    \end{itemize} 
\end{proposition}

We note that the application of the iterative step has two
formal conditions (the order of tangency to the identity, and
the Periodic Orbit Condition). There is also a  quantitative
condition \eqref{iterativeassumption}, telling us that the domain
loss cannot be too small with respect to the error.

To show that the iterative step can indeed be repeated, we will need
to recover the conditions. The formal conditions are recovered
automatically, but \eqref{iterativeassumption} and \eqref{est-C} will require
to specify the sequence of domain losses and show
that, under appropriate assumptions, the error decreases fast enough
so that \eqref{iterativeassumption} is maintained through the iteration.
This is very standard in KAM theory.

\begin{proof} 
Consider 
$\eta_{\eps}=\Id+ \sum_{j\geq L} \eta_{\eps}^{j}\eps^j$, 
and define $\phi_{\eps}$ ans $\tilde \eta_{\eps}$ as in Section \ref{s_iter_step_formal}. By definition, since $\eta_{\eps}$ satisfies the POC condition, $\tilde \eta_{\eps}$ does so as well.
The fact that $ \tilde \eta_\eps = \Id +  \tilde\eta_\eps^{[\ge 2L]}$ was proved in Remark \ref{r_eta_tilde}.  To prove estimate \eqref{est-C},  we will 
use formula  \eqref{eta-id-formal} and estimate the terms $A$, $B$ and $C$ in it.
\begin{itemize}
 
 \item By the Cauchy estimates recalled in Lemma \ref{lem:Cauchy}, 
 there exists a constant $c_1=c_1(d,\rho_0)$, such that
in formula \eqref{eta-id-formal}  we have:
$$
\begin{aligned}
\|  C \|_{\rho-\dt}= \| \eta_\eps^{[\ge 2L]}  \|_{\rho-\dt}  
\leq  c_1 \dt^{-1} \sigma  \exp{ (-   \delta L/\rho_0) } .
\end{aligned}
$$
In the same way,
\[
\| \eta^{[L, 2L -1]}_\eps  \|_{\rho-\dt}  \leq c_1 \dt^{-1} \sigma  .
\]

  \item Since  $\eta_\eps$ satisfies the POC condition, Proposition   \ref{POC->CPOC} implies that $\eta^{[L, 2L-1]}_\eps$ satisfies the CPOC condition. 
  By Proposition \ref{POC-FCequivalence}, condition CPOC implies condition FC for $\eta_\eps^{[L, 2L -1]}$.  
  Under the FC condition for $\eta_\eps^{[L, 2L -1]}$, Lemma~\ref{commutative_estimates} permits us to estimate:
  $$
  \|  \phi^{[L, 2L -1]}_\eps \|_{\rho-\dt}  \leq c_2 \dt^{-s} \| \eta^{[L, 2L-1]}_\eps  \|_{\rho-\dt}  \leq c_2 \dt^{-s} \si ,
 $$
 where $c_2$ only depends on the dimension.

\item  By  \eqref{iterativesol}, we can express  $ \phi_\eps^{[L, 2L -1]}\circ f =
     \eta_\eps^{[L, 2L -1]} + \phi_\eps^{[L, 2L -1]}$. This permits us to estimate:
\[
\begin{aligned}
\|\phi_\eps^{[L, 2L -1]} \circ f \|_{\rho -\delta}  \leq &  \|\eta_\eps^{[L, 2L -1]}  \|_{\rho -\delta} +\|\phi_\eps^{[L, 2L -1]} \|_{\rho -\delta}   \\
\leq & c_1 \dt^{-1} \sigma+ c_2 \delta^{-s} \sigma  \leq c_3 \dt^{-s} \si .
\end{aligned}
\]
Note that $c_3$ only depends on the dimension.  
 
 \item In the formulation of the proposition choose  $c'=c_3$. 
Assuming $c_3 \dt^{-s} \si\leq 1/100$, we can use Taylor's formula for $(1+x)^{-1}=1-x+O(x^2)$ applied to a matrix argument. Together with the boundedness of $\|\Id+\phi_\eps^{[L, 2L -1]} \|_{\rho -\delta} $ and $\|\eta_\eps \|_{\rho -\delta} $, this implies:
$$
\begin{aligned}
\|  A \|_{\rho- \dt}  \leq  & c_4  \| \left(\Id +  \phi^{[L, 2L -1]}_\eps\circ f \right)^{-1}  - \left( \Id  - \phi^{[L, 2L -1]}_\eps\circ f \right)  \|_{\rho- \dt} \\
\leq
& c_5 \| \phi^{[L, 2L -1]}_\eps\circ f  \|^2_{\rho- \dt} \leq c_6 \dt^{-2s} \si^2.
\end{aligned}
$$ 
\item Using the estimates for 
$\phi_\eps$ and $\eta_\eps$ above, we get
$$
\|  B \|_{\rho-\dt}  \leq c_6 \dt^{-2s} \si^2.
$$
\item  Now the upper bound of
  $\tilde \eta_\eps - \Id$ in the formula \eqref{eta-id-formal} is just a sum of the upper bounds above. We define $c$ to be the maximum of the relevant constants $c_1, \dots c_6$.
 
\end{itemize} 
\end{proof} 

\subsection{Convergence of the iterative procedure}
\label{sec:convergence}
In this section, we show that, starting from a sufficiently small
error, we can repeat on applying the iterative step and
the accumulated transformation converges.
The following calculation will be used in the construction.
\begin{lemma}\label{l_est_const} 
Let $a\geq 0$, $b>1$ and $c>0$ be fixed. Given $\ga_0$, 
define  a sequence of numbers $(\ga_n)_n$, $n\in\NN$, by a recursive formula 
 $$
 \ga_{n+1} = c \ga_n^2 \exp (a n) +c \ga_n\exp(-b^n).
 $$ 

For any choice of $(\la, p)$ such that $\la>0$ and $p\in(1,b]$, $p<2$,  
one can find $\Gamma_0=\Gamma_0(\la, p,a,b,c)$  so that for any $\ga_0\leq \Gamma_0$
the corresponding sequence satisfies:
\begin{equation}\label{gamma_exp_estimate}
\ga_n \leq \la  \exp (-p^{n}-an) \quad \text{ for all }  n\in \NN.
\end{equation} 
\end{lemma}

\noindent {\it Proof: }
Fix $(\la, p)$ as above. 
Since $p<2$, there exists an $N_0$ such that for all $n\geq N_0$
we have: 
$$
c(\la + 1)\exp(-2p^n) <  \exp(-p^{n+1} - a).
$$
Pick  $\Gamma_0=\Gamma_0(\la, p,a,b,c)$  so that  \eqref{gamma_exp_estimate} holds  for all $n\leq N_0$ and for all $\ga_0\leq \Gamma_0$.
This is possible since $\ga_0$ is a multiple in the recurrence relation.

The statement is proved by induction in $n$ with $n=N_0$ being the base.
Suppose that \eqref{gamma_exp_estimate} holds for some $n\geq N_0$.  Then  
\[
\begin{aligned}
 \ga_{n+1}& = c \ga_n^2 \exp (an) +c \ga_n \exp(-b^n) \\
 &\leq  c \la^2 \exp (-2p^{n} - 2an)+c \la  \exp(-p^n -b^n - an ) \\
 &\leq c( \la^2+\la)  \exp (-2p^{n} - an) \\
& \leq   \la \exp  (-p^{n+1} -a(n+1)).\\
 \end{aligned}
 \] 
Hence,   \eqref{gamma_exp_estimate} holds for $n+1$, and thus for all $n\in\NN$.   

\qed


\begin{remark}\label{rem_rescale}
We start with $\rho=\rho_0$,  $\eta_{\eps}\in \cA_\rho^\eps$ and $\si>0$  such  that 
$$
\| \eta_\eps -\Id \|_{\rho} =\si .
$$ 
To start the recursion, we need $\si$ to be small. This is reached by rescaling the parameter. 
Namely, for a real $\lambda > 0 $, consider the family
$\hat \eta_\eps  = \eta_{\lambda \eps}$. The analyticity domain of the rescaled family 
$\hat \eta_\eps$ in $(\th, I)$  is
the same as that of $\eta_\eps$; also, $\hat \eta_\eps$ satisfies
the Periodic Orbit Condition if and only if $\eta_\eps$ does.

At the same time, since $\eta_0 = \hat \eta_0  = \Id$,
we have
\[
\| \hat  \eta - \Id \|_\rho \le C \lambda^{-1}.
\]
Hence, by choosing $\lambda$ large enough, 
we can make $  \| \hat  \eta - \Id \|_\rho $ as small as desired.

Below we prove that, if
the initial error is small enough, the iterative method converges, and thus
 $\hat \eta_\eps$ is a coboundary with an analytic conjugacy in a domain defined by 
$\rho' < \rho$.
If $\hat \eta_\eps$ satisfies \eqref{coboundary}  with a conjugacy $\hat \phi_\eps$,
we see that $\eta_\eps$ is a coboundary with a conjugacy $\phi_\eps = \hat \phi_{\lambda^{-1} \eps}$.
This is enough to prove the conclusions of Theorem~\ref{main}.

Note that the domain of convergence of $\phi_\eps$ in the $\eps$ variable can be
significantly smaller than  that for $\eta_\eps$. 
\end{remark}

{\bf Basic notations and assumptions }
\begin{itemize}
\item Let $c$ be the constant from Proposition \ref{prop:iterative}, and let $s=d+1$.

\item Choose $ \dt_0  \leq   \rho_0 /8 $, and define $\delta_n = \dt_0 (3/4)^{n}$, 
$$
\rho_n= \rho_{0} - \sum_{j=1}^{n-1} \delta_{j-1}, \quad \rho_\infty= \rho_{0} - \sum_{j=1}^{\infty} \delta_{j-1}.
$$
Clearly, 
$$
\rho_\infty \geq \rho_{0}/2.
$$

\item Define $L_0$ so that $(\dt_0 L_0/\rho_0 ) \geq 1$, and let
$$
 L_n=L_0 2^n .
$$
Note that $(\dt_n L_n/\rho_0 ) = (\dt_0 L_0/\rho_0  )(3/2)^{n} \geq (3/2)^{n}$.

\item In Lemma \ref{l_est_const}, let  $a$ be such that 
$$ 
\exp(na)=\dt_n^{-2s} 
$$ 
(in this case $a$ is close to $2s \ln 4/3$), take 
$p=b=3/2$, $c$ from Proposition \ref{prop:iterative} and  $\la=1/100$.
Let $\Gamma_0(\la, p,a,b,c)$ be the constant given by Lemma \ref{l_est_const}. Fix  $\ga_0=\Gamma_0(\la, p,a,b,c)$, 
and let $(\gamma_n)$ be the corresponding sequence from Lemma  \ref{l_est_const}. Then \eqref{gamma_exp_estimate} reads as
$$
\gamma_n \leq \frac1{100} \delta^{2s}_n \exp (-(3/2)^n).
$$

\end{itemize} 


{\bf Proof of the theorem. }

Let us start the iterative procedure. 
We will drop the subscript $\eps$ in the notation of $\eta_\eps$ and $\phi_\eps$ for notational simplicity. We add a subindex $n$ to indicate the estimates on the $n$-th  iterative step. In particular, the starting function is $\eta_0$ satisfying POC condition,  and 
$$
\| \eta_0 -\Id \|_{\rho_0} =\si_0.
$$
By Remark \ref{rem_rescale} we can assume that  $\si_0\leq \ga_0$ defined above. 

At the $n$-th step, given $\eta_n$ satisfying POC, we construct $\phi_n$ by Proposition \ref{prop:iterative}. 
Then $\eta_{n+1}:=\phi_n\circ f\cdot \eta_n \cdot \phi_n$ satisfies POC condition.

Suppose that at each step estimate  \eqref{iterativeassumption} is satisfied. Note that estimate  \eqref{est-C} coinsides with the iterative relation defining the sequence $(\ga_n)$. This implies that for all $n$  we have:
$$
\si_n:=\| \eta_{n} - \Id \|_{\rho_n} \leq \ga_n.
$$
For our choice of $\si_0$,  Lemma  \ref{l_est_const} guarantees that 
$$
\si_n\leq\gamma_n\leq   \frac1{100} \dt_n^{2s} \exp (-(3/2)^{n}).
$$
This implies, in particular, that $c\si_n \dt_n^{-s} \leq  \frac1{100}$ for each $n$, i.e., \eqref{iterativeassumption} holds, and the inductive step can indeed be repeated.
 
Moreover, the same estimate, combined with   Proposition \ref{prop:iterative},   implies that, for an appropriate constant $c_1$, we have:
$$
\|\phi_n\|_{\rho_n-\dt_n} \leq c\dt_n^{-s} \si_n < \frac1{100}\dt_n^{s}< c_1 \left( \frac34\right)^{ns}.
$$

\medskip
 
 Let
 $$
 \Phi_n=\phi_0 \cdot \phi_{1}  \dots  \phi_n,
 $$
 in which case
 $$
 \Phi_n^{-1}\circ f=\phi_n^{-1}\circ f   \dots  \phi_{1}^{-1}\circ f \cdot  \phi_0^{-1}\circ f,
 $$
 and let
 $$
 \eta_n=\Phi_n^{-1}\circ f  \cdot \eta \cdot \Phi_n.
 $$
 Using formulas 
 $$
  \|  \phi_n -\Id  \|_{\rho_n-\dt_n}  =  \|  \phi^{[L_n, 2L_n -1]} \|_{\rho_n-\dt_n}  \leq c_1\left( \frac34\right)^{ns}
 $$ 
 and   
$$
  \|  (\phi_n)^{-1}\circ f -\Id  \|_{\rho_n -\dt_n}  \leq c_2 \left( \frac34\right)^{ns},
$$
we estimate:
$$
 \|  \Phi_n   \|_{\rho_0/2} \leq \prod_{j=1}^n   \|  \phi_j  \|_{\rho_0/2} \leq \exp(\sum_{j=1}^n   \|  \phi_j -\Id  \|_{\rho_0/2})
 \leq \exp(\sum_{j=1}^n \left( \frac34\right)^{ns}) .
$$
The latter is bounded  uniformly in $n$. By a similar argument, $ \|  \Phi_n ^{-1} \circ f   \|_{\rho_0/2} $ is bounded uniformly in $n$. 
We conclude that: $\Phi_n \in \cA_{\frac12 \rho_0}^\eps$,  $\Phi_n^{-1} \circ f \in \cA_{\frac12\rho_0 }^\eps$,
and the following limits exist:
$$ 
\phi^{\infty}:=\lim_{n\to\infty} \Phi_n \in \cA_{\frac12 \rho_0 }^\eps, \quad 
 (\phi^{\infty})^{-1}\circ f = \lim_{n\to\infty}(\Phi_n)^{-1}\circ f \in \cA_{\frac12 \rho_0 }^\eps.
$$

Moreover,
\[
\| (\Phi_n \circ f)^{-1}\cdot  \eta \cdot  \Phi_n - \Id \|_{\rho/2}\leq \si_n,
\]
which goes to zero when $n\to\infty$. Hence,
\[
(\phi^\infty\circ f)^{-1} \cdot \eta \cdot  \phi^\infty = \Id .
\]

\section{Adapting the proof to the case where $\G$ is a Lie group}
\label{sec:Lie}

In this section we explain (mostly notational) changes that permit us to adapt 
the above result to the case when $\G$ is a
Lie group (rather than an algebra of operators). An important example in
applications is
$\G=Sp(n, \complex)$.  

Since our results involve assumptions that the group elements are close to identity,
it is natural to use the exponential mapping which provides
an analytic diffeomorphism from a  neighbohood of  zero in the Lie algebra
to a nieghborhood of the identity in
the group. We write 
\begin{equation}\label{notation}
  \begin{split}
    & \eta_\eps  =  \exp(  \tilde \eta_\eps) \equiv
    \exp\left(\sum_{j = 1}^\infty {\tilde \eta}_j \eps^n  \right),\\
        & \phi_\eps  =  \exp(  \tilde \phi_\eps) \equiv
    \exp\left(\sum_{j = 1}^\infty {\tilde \phi}_j \right) \eps^n .
  \end{split}
\end{equation}
We note that since $\exp$ is an analytic diffeomorphism,
the family $\eta_\eps$ is analytic  if and
only if $\tilde \eta_\eps$ is convergent.
Similarly, $\tilde \phi_\eps$ is a formal
power series if and only if $\phi_\eps$ is.

We will refer to the elements of the Lie algebra corresponding to
the elements in the Lie group as the logarithms, and denote them
by $\tilde{}$.

As it is well known, the group multiplication in the Lie group becomes
(approximately) the sum  in the Lie algebra.
In our case the error is quadratic:
$$
\exp(a) \cdot \exp(b) = \exp(a+b) + O( (|a| + |b|)^2) .
$$
Indeed, there are asymptotic formulas for the error in terms of
commutators (Campbell-Hausdorff  formulas). 

The key observation  for our problem is that if $\tilde \eta_\eps= O(\eps^N)$ and
$\tilde \gamma_\eps  =  O(\eps^N)$, we have:
\begin{equation}\label{logproduct}
\eta_\eps \gamma_\eps =  \exp( \tilde \eta_\eps + \tilde \gamma_\eps + E ),
\end{equation} 
where $E = O(\eps^{2N})$ and its analytic norm can be bounded
by the sums of the squares of the norms of $\eta_\eps$ and $\gamma_\eps$.

With these observations, the proof
of an analogue of Theorem~\ref{main} goes through with only minimal changes.
We inductively consider $\tilde \eta_\eps$ vanishing to order $L$ in $\eps$, 
such that $\eta_\eps $ satisfies the Periodic Orbit Condition.
Matching powers in the product, we again obtain that
$\tilde \eta^{[L,2L-1]}_\eps$ satisfies the Commutative Periodic Orbit Condition.  
Hence, we can determine $\tilde \phi^{[L,2L-1]}_\eps$ solving
the commutative cohomology equation.

We see, proceeding as in the proof of Theorem ~\ref{main},
that the new $\eta_\eps$ is quadratically small with respect to the one at the previous step. 
The terms that appear in the estimates 
are the same as those considered there, supplemented by a
few terms coming from the $E$ in \eqref{logproduct}. Since the terms coming
from \eqref{logproduct}
 satisfy uniform quadratic estimates, they  do not affect the estimates 
 of the inductive step, and the convergence of the Nash-Moser method remains exactly the same.


\begin{thebibliography}{dlLMM86}

\bibitem[dlL96]{Llave96}
Rafael de~la Llave.
\newblock On necessary and sufficient conditions for uniform integrability of
  families of {H}amiltonian systems.
\newblock In {\em International {C}onference on {D}ynamical {S}ystems
  ({M}ontevideo, 1995)}, volume 362 of {\em Pitman Res. Notes Math. Ser.},
  pages 76--109. Longman, Harlow, 1996.

\bibitem[dlL97]{Llave97}
R.~de~la Llave.
\newblock Analytic regularity of solutions of {L}ivsic's cohomology equation
  and some applications to analytic conjugacy of hyperbolic dynamical systems.
\newblock {\em Ergodic Theory Dynam. Systems}, 17(3):649--662, 1997.

\bibitem[dlLMM86]{LlaveMM86}
R.~de~la Llave, J.~M. Marco, and R.~Moriy\'{o}n.
\newblock Canonical perturbation theory of {A}nosov systems and regularity
  results for the {L}iv\v{s}ic cohomology equation.
\newblock {\em Ann. of Math. (2)}, 123(3):537--611, 1986.

\bibitem[dlLS21]{dlLS}
Rafael de~la Llave and Maria Saprykina.
\newblock Convergence of the birkhoff normal form sometimes implies convergence
  of a normalizing transformation.
\newblock {\em Ergodic Theory and Dynamical Systems}, page 1–22, 2021.

\bibitem[dlLW10]{LlaveW10}
Rafael de~la Llave and Alistair Windsor.
\newblock Liv\v{s}ic theorems for non-commutative groups including
  diffeomorphism groups and results on the existence of conformal structures
  for {A}nosov systems.
\newblock {\em Ergodic Theory Dynam. Systems}, 30(4):1055--1100, 2010.

\bibitem[Kal11]{Kalinin11}
Boris Kalinin.
\newblock Liv\v{s}ic theorem for matrix cocycles.
\newblock {\em Ann. of Math. (2)}, 173(2):1025--1042, 2011.

\bibitem[KN11]{KatokN11}
Anatole Katok and Viorel Ni\c{t}ic\u{a}.
\newblock {\em Rigidity in higher rank abelian group actions. {V}olume {I}},
  volume 185 of {\em Cambridge Tracts in Mathematics}.
\newblock Cambridge University Press, Cambridge, 2011.
\newblock Introduction and cocycle problem.

\bibitem[Liv71]{Livsic71}
A.~N. Liv\v{s}ic.
\newblock Certain properties of the homology of {$Y$}-systems.
\newblock {\em Mat. Zametki}, 10:555--564, 1971.

\bibitem[Liv72]{Livsic72}
A.~N. Liv\v{s}ic.
\newblock Cohomology of dynamical systems.
\newblock {\em Izv. Akad. Nauk SSSR Ser. Mat.}, 36:1296--1320, 1972.

\bibitem[Mos67]{Moser67}
J\"{u}rgen Moser.
\newblock Convergent series expansions for quasi-periodic motions.
\newblock {\em Math. Ann.}, 169:136--176, 1967.

\bibitem[NT95]{NiticaT95}
Viorel Ni\c{t}ic\u{a} and Andrei T\"{o}r\"{o}k.
\newblock Cohomology of dynamical systems and rigidity of partially hyperbolic
  actions of higher-rank lattices.
\newblock {\em Duke Math. J.}, 79(3):751--810, 1995.

\bibitem[NT01]{NiticaT01}
Viorel Ni\c{t}ic\u{a} and Andrei T\"{o}r\"{o}k.
\newblock Local rigidity of certain partially hyperbolic actions of product
  type.
\newblock {\em Ergodic Theory Dynam. Systems}, 21(4):1213--1237, 2001.

\bibitem[NT02]{NiticaT02}
Viorel Ni\c{t}ic\u{a} and Andrei T\"{o}r\"{o}k.
\newblock On the cohomology of {A}nosov actions.
\newblock In {\em Rigidity in dynamics and geometry ({C}ambridge, 2000)}, pages
  345--361. Springer, Berlin, 2002.

\bibitem[Poi99]{Poincare99}
H.~Poincar{\'e}.
\newblock {\em Les m\'ethodes nouvelles de la m\'ecanique c\'eleste}, volume 1,
  2, 3.
\newblock Gauthier-Villars, Paris, 1892--1899.

\bibitem[{Vee}86]{Veech86}
William~A. {Veech}.
\newblock {Periodic points and invariant pseudomeasures for toral
  endomorphisms}.
\newblock {\em {Ergodic Theory Dyn. Syst.}}, 6:449--473, 1986.

\bibitem[Wil13]{Wilkinson13}
Amie Wilkinson.
\newblock The cohomological equation for partially hyperbolic diffeomorphisms.
\newblock {\em Ast\'{e}risque}, (358):75--165, 2013.

\end{thebibliography}

\end{document}